\documentclass[12pt]{article}
\usepackage{latexsym,amsfonts,amssymb}
\usepackage[dvips]{color}
\usepackage{colordvi,multicol}

\def \cal{\mathcal}

\newtheorem{thm}{Theorem}[section]
\newtheorem{cor}[thm]{Corollary}
\newtheorem{lem}[thm]{Lemma}
\newtheorem{pro}[thm]{Proposition}
\newtheorem{defi}[thm]{Definition}
\newtheorem{rem}[thm]{Remark}

\date{}
\begin{document}
\title{\bf  Two Theorems on Hunt's Hypothesis  (H) for  Markov Processes}
\author{}

\maketitle

 \centerline{Ze-Chun Hu} \centerline{\small College
of Mathematics, Sichuan University} \centerline{\small  Chengdu,
610064, China} \centerline{\small E-mail: zchu@scu.edu.cn}

\vskip 0.5cm \centerline{Wei Sun} \centerline{\small Department of
Mathematics and Statistics, Concordia University}
\centerline{\small Montreal, H3G 1M8, Canada} \centerline{\small
E-mail: wei.sun@concordia.ca}

\vskip 0.5cm \centerline{Li-Fei Wang} \centerline{\small College
of Mathematics and Information Science, Hebei Normal University}
\centerline{\small Shijiazhuang, 050024, China} \centerline{\small
E-mail: flywit1986@163.com}

\begin{abstract}

\noindent Hunt's hypothesis (H) and the related Getoor's
conjecture is one of the most important problems in the basic
theory of Markov processes.  In this paper, we investigate the
invariance of Hunt's hypothesis (H) for Markov processes under two
classes of transformations, which are change of measure and
subordination. Our first theorem shows that for two standard
processes $(X_t)$ and $(Y_t)$, if $(X_t)$ satisfies (H) and
$(Y_t)$ is locally absolutely continuous with respect to $(X_t)$,
then $(Y_t)$ satisfies (H). Our second theorem shows that a
standard process $(X_t)$ satisfies (H) if and only if
$(X_{\tau_t})$ satisfies (H) for some (and hence any) subordinator
$(\tau_t)$ which is independent of $(X_t)$ and has a positive
drift coefficient. Applications of the two theorems are given.
\end{abstract}



%

\section{Introduction and main results}

Markov processes are a very important class of stochastic
processes and have broad applications in many areas. Although
their basic theory has been well-developed, there are still some
challenging unsolved problems. One of them is Hunt's hypothesis
(H) and the related Getoor's conjecture. This hypothesis plays a
crucial role in the potential theory of Markov processes.

The study of the relationship between Markov processes and
potential theory has a long history. It dated back to the work of
Feller, Kac, Kakutani, et al. and Doob's extensive study of the
relationship between Brownian motion and classical potential
theory \cite{Doob3}. We refer the reader to the monographs of Port
and Stone \cite{PS78} and Doob \cite{Do84} for full descriptions.
In 1957 and 1958, Hunt laid the foundation for probabilistic
potential theory in his celebrated papers \cite{Hu57a,Hu57b,Hu58}.
His work revealed the deep connection between general Markov
processes and potential theory. We refer the reader to Blumenthal
and Getoor \cite{BG68}, Aikawa and Ess\'{e}n \cite{AE96}, and
Beznea and Boboc \cite{BB04} for more introduction to Markov
processes and potential theory.

Let $E$ be a locally compact space with a countable base (LCCB)
and $X=(X_t, P^x)$ be a standard Markov process on $E$ as
described in Blumenthal and Getoor \cite{BG68}. Denote by ${\cal
B}$ and ${\cal B}^n$ the family of all Borel measurable subsets
and nearly Borel measurable subsets of $E$, respectively. For
$D\subset E$, we define the first hitting time of $D$ by
$$
T_D=\inf\{t>0:X_t\in D\}.
$$

A set $D\subset E$ is called thin if there exists a set $C\in
{\cal B}^n$ such that $D\subset C$ and {$P^x(T_C=0)=0$} for any
$x\in E$. $D$ is called semipolar if
$D\subset\bigcup_{n=1}^{\infty}D_n$ for some thin sets
$\{D_n\}_{n=1}^{\infty}$. $D$ is called polar if there exists a
set $C\in {\cal B}^n$ such that $D\subset C$ and
$P^x(T_C<\infty)=0$ for any $x\in E$. Let $m$ be a measure on
$(E,{\cal B})$. $D$ is called $m$-essentially polar if there
exists a set $C\in {\cal B}^n$ such that $D\subset C$ and
$P^m(T_C<\infty)=0$. Hereafter
$P^m(\cdot):=\int_EP^x(\cdot)m(dx)$.

Hunt's hypothesis (H) says that ``every semipolar set of $X$ is
polar". It has various equivalent expressions in probabilistic
potential theory. For example, it is known that if $X$ is in
duality with another standard process $\hat{X}$ on $E$ with
respect to a $\sigma$-finite reference measure $m$, then (H) is
equivalent to many potential principles for Markov processes,
e.g., the bounded maximum principle (cf. \cite{Fi90}); (H) holds
if and only if the fine and cofine topologies differ by polar sets
(\cite[Proposition 4.1]{BG70} and \cite[Theorem 2.2]{G83}); (H)
holds if and only if every natural additive functional of $X$ is
in fact a continuous additive functional (\cite[Chapter
VI]{BG68}).

In spite of its importance, (H) has been verified only in special
situations. Some fifty years ago, Getoor conjectured that
essentially all L\'{e}vy processes  satisfy (H), except for some
extremely nonsymmetric cases like uniform motions. This conjecture
stills remains open and is a major unsolved problem in the
potential theory for L\'{e}vy processes. The reader is referred to
Kesten \cite{Ke69}, Port and Stone \cite{PS69}, Blumenthal and
Getoor \cite{BG70}, Bretagnolle \cite{Br71}, Forst \cite{F75},
Kanda \cite{Ka76}, Rao \cite{R77}, Kanda \cite{Ka78}, Glover and
Rao \cite{GR86}, and Rao \cite{R88} for the results that were
obtained before 1990. The reader is also referred to Hu and Sun
\cite{HS12}, Hu et al. \cite{HSZ15}, Hu and Sun \cite{HS16}, and
Hu and Sun \cite{HS18} for the recent results on Getoor's
conjecture.

We would like to call the reader's attention to the related work
on Hunt's hypothesis (H) and the potential theory of Markov
processes. These include Silverstein \cite{Si77}, Hawkes
\cite{Ha79}, Glover and Rao \cite{GR88}, Fitzsimmons and Kanda
\cite{FK92},   Fitzsimmons \cite{Fi01}, Beznea et al.
\cite{BCR11}, Han et al. \cite{HMS11}, Fitzsimmons \cite{Fi14},
Hansen and Netuka \cite{HN16}.

In this paper, we study Hunt's hypothesis (H) for general
Markov processes from the point of view of transformations. We
will present two new theorems on (H), which imply that various
classes of Markov processes satisfy (H).

We fix an isolated point
$\Delta$ which is not in $E$ and write $E_{\Delta}=E\cup
\{\Delta\}$. Consider the following objects:

(i) $\Omega$ is a set and $\omega_{\Delta}$ is a distinguished point of $\Omega$.

(ii) For $0\le t\le\infty$, $Z_t:\Omega \rightarrow E_{\Delta}$ is
a map such that if $Z_t(\omega)=\Delta$  then $Z_s(\omega)=\Delta$
for all $s\ge t$, $Z_{\infty}(\omega)=\Delta$ for all
$\omega\in\Omega$, and $Z_0(\omega_{\Delta})=\Delta$.

(iii) For $0\le t\le\infty$, $\theta_t:\Omega\rightarrow\Omega$ is a map such that $Z_s\circ \theta_t=Z_{s+t}$ for all $s,t\in [0,\infty]$, and $\theta_{\infty}\omega=\omega_{\Delta}$ for all $\omega\in\Omega$.

We define in $\Omega$ the $\sigma$-algebras
$\mathcal{F}^0=\sigma(Z_t: t\in [0,\infty])$ and
$\mathcal{F}^0_t=\sigma(Z_s: s\leq t)$ for $0\le t<\infty$. Denote
$$
\zeta(\omega)=\inf\{\omega:Z_t(\omega)=\Delta\},\ \ \omega\in \Omega.
$$
Let $m$ be a measure on $(E,{\cal B})$. We define
$$(H_m):\ \mbox{every semipolar set is}\ m\mbox{-essentially\ polar}.$$
Note that if a standard process $X$ has resolvent densities with
respect to $m$, then $X$ satisfies (H) if and only if $X$
satisfies $(H_m)$ (cf. \cite[Propositions II.2.8 and
II.3.2]{BG68}).

Now we can state the first main result of this paper.
\begin{thm}\label{thm-2}
Let $X=(\Omega,\mathcal{M}^X,\mathcal{M}^X_t,Z_t,\theta_t,P^x)$ and
$Y=(\Omega,\mathcal{M}^Y,\mathcal{M}^Y_t,Z_t,\theta_t,Q^x)$ be two
standard processes on $E$ such that $\mathcal{M}^X\cap \mathcal{M}^Y\supset \mathcal{F}^0$ and $\mathcal{M}^X_t\cap \mathcal{M}^Y_t\supset \mathcal{F}_t^0$ for $0\le t<\infty$.

(i) Suppose that $X$ satisfies (H) and for any $x\in E$ and $t>0$, $Q^x|_{\cal{F}_t^0}$ is absolutely
continuous with respect to $P^x|_{\cal{F}_t^0}$ on $\{t<\zeta\}$. Then $Y$ satisfies (H).

(ii) Suppose that $X$ satisfies ($H_m$) for some measure $m$ on $(E,{\cal B})$ and for any $x\in E$ and $t>0$, $Q^x|_{\cal{F}_t^0}$ is absolutely
continuous with respect to $P^x|_{\cal{F}_t^0}$ on $\{t<\zeta\}$. Then $Y$ satisfies ($H_m$).
\end{thm}

A subordinator $\tau=(\tau_t)$ is a 1-dimensional increasing
L\'{e}vy process with $\tau_0=0$. Let $X=(X_t)$ be a standard
process on $E$ and $\tau$ be a subordinator  which is independent
of $X$. The standard process $(X_{\tau_t})$ is called the
subordinated process of $(X_t)$. The idea of subordination
originated from Bochner (cf. \cite{Bo55}). Our second theorem is
motivated by the following remarkable result of Glover and Rao.
\begin{thm} \label{thm-3}  (Glover and Rao \cite{GR86}) Let $(X_t)$
be a standard process on $E$ and $(\tau_t)$ be a subordinator
which is independent of $X$ and satisfies (H). Then $(X_{\tau_t})$
satisfies (H).
\end{thm}

It is known that if a subordinator $(\tau_t)$ satisfies (H), then
it must be a pure jump subordinator, i.e., its drift coefficient
equals 0 (\cite[Proposition 1.6]{HS12}). Up to now, it is still
unknown if any pure jump subordinator satisfies (H). We present
the following new theorem on the equivalence between (H) for $X$
and (H) for its time changed process.

\begin{thm}\label{thm-1}
Let $(X_t)$ be a standard process on $E$ and $m$ be a measure on
$(E,{\cal B})$. Then,

(i) $(X_t)$ satisfies (H) if and only if $(X_{\tau_t})$ satisfies (H)
for some (and hence any) subordinator $(\tau_t)$ which is
independent of $(X_t)$ and has a positive drift coefficient.

(ii) $(X_t)$ satisfies ($H_m$) if and only if $(X_{\tau_t})$ satisfies ($H_m$)
for some (and hence any) subordinator $(\tau_t)$ which is
independent of $(X_t)$ and has a positive drift coefficient.
\end{thm}

The rest of this paper is organized as follows. In Sections
\ref{sec3} and \ref{sec2}, we give the proofs of Theorems
\ref{thm-2} and \ref{thm-1}, respectively. Applications of
Theorems \ref{thm-2} and \ref{thm-1} will be given in Sections
\ref{sec4} and \ref{sec5}, respectively.

\section{Proof of Theorem \ref{thm-2}}\label{sec3}\setcounter{equation}{0}

\subsection{Preliminary lemmas}

Before proving Theorem \ref{thm-2}, we give two lemmas. The first one is a well-known result, but we present it here for the reader's convenience.

\begin{lem}\label{rlem-3.5}
Let $X=(X_t,P^x)$ be a standard process on $E$ and $m$ be a
measure on $(E,{\cal B})$.

(i) If any thin set $A\in {\cal B}$ is polar, then $X$ satisfies
(H).

(ii) If any thin set $A\in {\cal B}$ is $m$-essentially polar,
then $X$ satisfies ($H_m$).
\end{lem}
{\bf Proof.} The proof of (ii) is similar to that of (i). So we only
prove (i) below.

Let $A$ be a semipolar set. We will show that $A$ is a polar set.
By the definitions of semipolar set and polar set, we may assume
without loss of generality that $A$ is a thin nearly Borel
measurable set. Then, we have $P^x(T_A=0)=0$ for any $x\in E$.

We fix an $x\in E$. By the definition of nearly Borel  measurable set, there exist two Borel sets $A_1$ and $A_2$ such that $A_1\subset A\subset A_2$ and
\begin{eqnarray*}
P^x(\{\omega: X_t(\omega)\in A_2-A_1\ \mbox{for some}\ t\in [0,\infty)\})=0,
\end{eqnarray*}
which implies that
\begin{eqnarray}\label{rlem-3.5-b}
T_{A_1}=T_A=T_{A_2},\ \ P^x\mbox{-a.s.}.
\end{eqnarray}
Then, $A_1\in {\cal B}$ is a thin set. By the assumption of the
lemma, we know that $A_1$ is a polar set. Therefore,
$P^x(T_{A_1}<\infty)=0$,  which together with (\ref{rlem-3.5-b})
implies that
$$
P^x(T_A<\infty)=0.
$$
Since $x\in E$ is arbitrary, this implies that $A$ is a polar
set.\hfill\fbox

Let $X=(\Omega,\mathcal{M}^X,\mathcal{M}^X_t,Z_t,\theta_t,P^x)$
and $Y=(\Omega,\mathcal{M}^Y,\mathcal{M}^Y_t,Z_t,\theta_t,Q^x)$ be
two standard processes on $E$ such that $\mathcal{M}^X\cap
\mathcal{M}^Y\supset \mathcal{F}^0$ and $\mathcal{M}^X_t\cap
\mathcal{M}^Y_t\supset \mathcal{F}_t^0$  for $0\le t<\infty$.
Suppose that for any $x\in E$ and $t>0$, $Q^x|_{\cal{F}_t^0}$ is
absolutely continuous with respect to $P^x|_{\cal{F}_t^0}$ on
$\{t<\zeta\}$.

Let $A$ be a subset of $E$. For $w\in \Omega$, we define
\begin{eqnarray*}
D_A(\omega):=\inf\{t\geq 0: Z_t(\omega)\in A\},\ \ T_A(\omega):=\inf\{t> 0: Z_t(\omega)\in A\}.
\end{eqnarray*}
Then,
\begin{eqnarray}\label{3.1-a}
\left\{\frac{1}{l}+D_A\left(\theta_{\frac{1}{l}}w\right)\right\}\downarrow
T_A(w)\ \ {\rm as}\ l\uparrow\infty.
\end{eqnarray}
By (\ref{3.1-a}), we know that for $t> 0$,
\begin{eqnarray}\label{lem-3.3-b-0}
\{T_A<t\}=\bigcup_{l=1}^{\infty}\left\{\frac{1}{l}+D_A\circ \theta_{\frac{1}{l}}<t\right\},
\end{eqnarray}
and
\begin{eqnarray}\label{lem-3.3-g-1}
\{T_A\geq t\}=\bigcap_{l=1}^{\infty}\left\{\frac{1}{l}+D_A\circ \theta_{\frac{1}{l}}\geq t\right\}.
\end{eqnarray}

\begin{lem}\label{lem-3.3}
Let $A\in {\cal B}$. Then, for any $x\in E$ and $t>0$,

(i) if $P^x(T_A< t,\ t<\zeta)=0$, then $Q^x(T_A< t,\ t<\zeta)=0$;

(ii) if $P^x(T_A>  t,\ t<\zeta)=0$, then $Q^x(T_A>  t,\ t<\zeta)=0$.
\end{lem}

We would like to point out that the proof of Lemma \ref{lem-3.3}
is far from trivial. Define $\mathcal{F}^X$ to be the completion
of $\mathcal{F}^0$ with respect to $\{P^{\mu}: \mu\ \mbox{is a
finite measure on}\ E_{\Delta}\}$, and define $\mathcal{F}^X_t$ to
be the completion of $\mathcal{F}^0_t$ in $\mathcal{F}^X$ with
respect to $\{P^{\mu}: \mu\ \mbox{is a finite measure on}\
E_{\Delta}\}$ for $t\in [0,\infty)$. Let $A\in {\cal B}$. By \cite[Theorem
I.10.7]{BG68}, we know that $T_A$ is a stopping time relative to
$\{\mathcal{F}_t^{X}\}$. Since $\cal{F}_t^X$ contains all null
sets in $\cal{F}^X$, in general, the  assumption that
$Q^x|_{\cal{F}_t^0}$ is absolutely continuous with respect to
$P^x|_{\cal{F}_t^0}$ on $\{t<\zeta\}$ does not imply that
$Q^x|_{\cal{F}_t^X}$ is absolutely continuous with respect to
$P^x|_{\cal{F}_t^X}$ on $\{t<\zeta\}$.

\noindent {\bf Proof of Lemma \ref{lem-3.3}.} {\it Step 1.}  If $A$ is an open set, then $T_A$ is an $\{\mathcal{F}^0_{t+}\}$ stopping time, where $\cal{F}^0_{t+}=\bigcap_{s>t}\cal{F}_s^0$. Note that
$$
\{T_A< t\}\bigcap\{t<\zeta\}\in{\cal F}^0_t\bigcap\{t<\zeta\},
$$
and
$$
\{T_A> t\}\bigcap\{t<\zeta\}=\bigcup_{l=1}^{\infty}\left(\{T_A> t\}\bigcap\left\{t+\frac{1}{l}<\zeta\right\}\right)\in\bigcup_{l=1}^{\infty}\left({\cal F}^0_{t+\frac{1}{l}}\bigcap\left\{t+\frac{1}{l}<\zeta\right\}\right).
$$
Thus, (i) and (ii) hold by the assumption that  $Q^x|_{\cal{F}_t^0}$ is absolutely
continuous with respect to $P^x|_{\cal{F}_t^0}$ on $\{t<\zeta\}$ for any $t>0$.

{\it Step 2.} Suppose that $A$ is a compact set. Then, there
exists a sequence $\{A_n\}$ of open sets satisfying $A_n\supset
\overline{A_{n+1}}$ and $\bigcap_{n\geq 1}A_n=A$. For $s>0$,
$\{s+D_{A_n}\circ \theta_s\}$ is an increasing sequence of
$\{\mathcal{F}^0_{t+}\}$ stopping times. Define
\begin{equation}\label{r11} D_s=\lim_{n\to\infty}(s+D_{A_n}\circ
\theta_s).
\end{equation}
 Then, $D_s$ is an $\{\mathcal{F}^0_{t+}\}$ stopping time.

Obviously, $D_s\leq s+D_A\circ \theta_s$. By the right continuity of $X$, we know that $X(s+D_{A_n}\circ \theta_s)\in \overline{A_n}$. Then, we obtain by the quasi-left continuity of $X$ that for any $x\in E$,
\begin{eqnarray*}
P^x(Z(D_s)\in A, D_s<\zeta)=P^x\left(\lim_{n\to\infty}Z(s+D_{A_n}\circ \theta_s)\in
\bigcap_{n\geq 1}\overline{A_n}=A, D_s<\zeta\right).
\end{eqnarray*}
It follows that $P^x(s+D_A\circ \theta_s\leq D_s, D_s<\zeta)=P_x(D_s<\zeta)$. Hence
\begin{eqnarray}\label{lem-3.3-a}
P^x(D_s\neq s+D_A\circ \theta_s)=0.
\end{eqnarray}
Similarly, we can show that for any $x\in E$,
\begin{eqnarray}\label{lem-3.3-b}
Q^x(D_s\neq s+D_A\circ \theta_s)=0.
\end{eqnarray}

If $P^x(T_A<t,\ t<\zeta)=0$, then we obtain by (\ref{lem-3.3-b-0}) that for any $l\in \mathbb{N}$,
\begin{eqnarray}\label{lem-3.3-c}
P^x\left(\frac{1}{l}+D_A\circ \theta_{\frac{1}{l}}<t,\ t<\zeta\right)=0.
\end{eqnarray}
For $l\in \mathbb{N}$, we have
\begin{eqnarray}\label{lem-3.3-d}
\left\{\frac{1}{l}+D_A\circ \theta_{\frac{1}{l}}< t\right\}&=&\left(\left\{D_{\frac{1}{l}}< t\right\}
\bigcap \left\{\frac{1}{l}+D_A\circ \theta_{\frac{1}{l}}=D_{\frac{1}{l}}\right\}\right)\nonumber\\
&&\bigcup \left(\left\{\frac{1}{l}+D_A\circ \theta_{\frac{1}{l}}< t\right\}\bigcap \left\{\frac{1}{l}+D_A\circ \theta_{\frac{1}{l}}\neq D_{\frac{1}{l}}\right\}\right).
\end{eqnarray}
Then, we obtain by (\ref{lem-3.3-a}), (\ref{lem-3.3-c}) and
(\ref{lem-3.3-d}) that
\begin{eqnarray}\label{lem-3.3-e}
P^x\left(D_{\frac{1}{l}}< t,\ t<\zeta\right)=0.
\end{eqnarray}

Note that $\{D_{\frac{1}{l}}< t\}\in \mathcal{F}_t^0$. By (\ref{lem-3.3-e}) and  the assumption that $Q^x|_{\cal{F}_t^0}$ is absolutely
continuous with respect to $P^x|_{\cal{F}_t^0}$ on $\{t<\zeta\}$, we get
\begin{eqnarray}\label{lem-3.3-f}
Q^x\left(D_{\frac{1}{l}}< t,\ t<\zeta\right)=0.
\end{eqnarray}
Then, we obtain by (\ref{lem-3.3-b}), (\ref{lem-3.3-d}) and (\ref{lem-3.3-f}) that for any $l\in \mathbb{N}$,
$$
Q^x\left(\frac{1}{l}+D_A\circ \theta_{\frac{1}{l}}<t,\ t<\zeta\right)=0,
$$
which together with (\ref{lem-3.3-b-0}) implies that
$$
Q^x(T_A<t,\ t<\zeta)=0.
$$
Therefore, (i) holds.

Now we show that (ii) holds. Since
$\{T_A>t,t<\zeta\}=\bigcup_{n=1}^{\infty}\{T_A\geq
t+\frac{1}{n},t+\frac{1}{n}<\zeta\}$, it is  sufficient to show
that for any $t>0$, $P^x(T_A\geq t,\ t<\zeta)=0$ implies that
$Q^x(T_A\geq  t,\ t<\zeta)=0$.

Suppose that $P^x(T_A\geq t,\ t<\zeta)=0$. Then, we obtain by
(\ref{lem-3.3-g-1}) that
\begin{eqnarray}\label{lem-3.3-h}
 P^x\left(\frac{1}{l}+D_A\circ \theta_{\frac{1}{l}}\geq
t,\ t<\zeta\right)\downarrow0\ \ {\rm as}\ l\uparrow\infty.
\end{eqnarray}
For $l\in \mathbb{N}$, we have
\begin{eqnarray}\label{lem-3.3-h-1}
\left\{\frac{1}{l}+D_A\circ \theta_{\frac{1}{l}}\geq  t\right\}&=&\left(\left\{D_{\frac{1}{l}}\geq  t\right\}\bigcap \left\{\frac{1}{l}+D_A\circ \theta_{\frac{1}{l}}=D_{\frac{1}{l}}\right\}\right)\nonumber\\
&&\bigcup \left(\left\{\frac{1}{l}+D_A\circ
\theta_{\frac{1}{l}}\geq  t\right\}\bigcap
\left\{\frac{1}{l}+D_A\circ \theta_{\frac{1}{l}}\neq
D_{\frac{1}{l}}\right\}\right).\ \
\end{eqnarray}
By (\ref{lem-3.3-a}), (\ref{lem-3.3-b}) and (\ref{lem-3.3-h-1}), we get
\begin{eqnarray}\label{lem-3.3-i}
P^x\left(\frac{1}{l}+D_A\circ \theta_{\frac{1}{l}}\geq t,\ t<\zeta\right)=P^x\left(D_{\frac{1}{l}}\geq t,\ t<\zeta\right),
\end{eqnarray}
and
\begin{eqnarray}\label{lem-3.3-j}
Q^x\left(\frac{1}{l}+D_A\circ \theta_{\frac{1}{l}}\geq t,\ t<\zeta\right)=Q^x\left(D_{\frac{1}{l}}\geq t,\ t<\zeta\right).
\end{eqnarray}

Note that $\{D_{\frac{1}{l}}\geq t\}\in \mathcal{F}_t^0$. Then, we obtain  by (\ref{lem-3.3-h}), (\ref{lem-3.3-i}), (\ref{lem-3.3-j}) and the assumption $Q^x|_{\cal{F}_t^0}$ is absolutely
continuous with respect to $P^x|_{\cal{F}_t^0}$ on $\{t<\zeta\}$ that
\begin{eqnarray*}
 Q^x\left(\frac{1}{l}+D_A\circ \theta_{\frac{1}{l}}\geq
t,\ t<\zeta\right)\downarrow0\ \ {\rm as}\ l\uparrow\infty,
\end{eqnarray*}
which together with (\ref{lem-3.3-g-1}) implies that
$Q^x(T_A\geq t,\ t<\zeta)=0$. Therefore, (ii) holds.

{\it Step 3.} Suppose that $A\in {\cal B}$ and $x\in E$. Let $t>0$ and $s\in (0,t)$. For $B\subset E$, define
$$
\wedge_{s,t}(B)=\{\omega: Z_u(\omega)\in B\ \mbox{for some}\ u\in [s,t]\}.
$$
Then,
\begin{equation}\label{import}
\{s+D_B\circ \theta_s<t\}\subset\wedge_{s,t}(B)\subset\{s+D_B\circ \theta_s\le t\},
\end{equation}
and
\begin{equation}\label{import2}
\{s+D_B\circ \theta_s\ge t\}\supset(\wedge_{s,t}(B))^c\supset\{s+D_B\circ \theta_s> t\}.
\end{equation}

Denote by $\mathcal{O}$ the family of all open subsets of $E$. If
$B\in \mathcal{O}$,  then $\wedge_{s,t}(B)=\bigcup_{r\in
\mathbb{Q}\cap (s,t)}\{X_r\in B\}\bigcup \{X_s\in B\}\bigcup
\{X_t\in B\}\in \mathcal{F}^0_t$. Define the set function $I$ on
$\mathcal{O}$ by
$$
I(B)=(P^x+Q^x)(\wedge_{s,t}(B)),\ B\in \mathcal{O}.
$$
Following \cite[Lemma A.2.6]{FOT94}, we can prove the following
proposition.

\begin{pro}\label{r5}
The set function $I$ on $\mathcal{O}$ satisfies the following conditions:

(I.1) $B_1,B_2\in \mathcal{O},B_1\subset B_2 \Rightarrow I(B_1)\leq I(B_2)$.

(I.2) $I(B_1\cup B_2)+I(B_1\cap B_2)\leq I(B_1)+I(B_2),\ B_1,B_2\in\mathcal{O}$.

(I.3) $B_n\in\mathcal{O},B_n\uparrow B\Rightarrow B\in \mathcal{O},I(B)=\lim_{n\to\infty}I(B_n)$.
\end{pro}

Define
$$
I^*(B)=\inf_{G\in \mathcal{O}: B\subset G}I(G),\ \ B\subset E.
 $$
By proposition \ref{r5} and \cite[Theorem A.1.2]{FOT94}, we know
that $I^*$ is a Choquet capacity. By (\ref{r11}),
(\ref{lem-3.3-a}) and (\ref{lem-3.3-b}), we find  that
$I^*(B)=(P^x+Q^x)(\wedge_{s,t}(B))$ if $B$ is a compact subset of
$E$. Then, we obtain by the Choquet theorem (cf. \cite[Theorem
A.1.1]{FOT94}) that there exist a decreasing sequence $\{A_n\}$ of
open sets and an increasing sequence $\{B_n\}$ of compact sets
such that $B_n\subset A\subset A_n$ and
$\lim_{n\to\infty}I(A_n)=\lim_{n\to\infty}I^*(B_n)$. Consequently,
we have
$$
\bigcup_{n=1}^{\infty}\wedge_{s,t}(B_n)\subset \wedge_{s,t}(A)\subset \bigcap_{n=1}^{\infty}\wedge_{s,t}(A_n)
$$
and
$$
(P^x+Q^x)\left(\bigcap_{n=1}^{\infty}\wedge_{s,t}(A_n)-
\bigcup_{n=1}^{\infty}\wedge_{s,t}(B_n)\right)=0,
$$
which implies that
\begin{eqnarray}\label{r22}
&&P^x(\wedge_{s,t}(A),\ t<\zeta)=P^x\left(\bigcap_{n=1}^{\infty}\wedge_{s,t}(A_n),\ t<\zeta\right),\nonumber\\
 &&Q^x(\wedge_{s,t}(A),\ t<\zeta)=Q^x\left(\bigcap_{n=1}^{\infty}\wedge_{s,t}(A_n),\ t<\zeta\right),
\end{eqnarray}
and
\begin{eqnarray}\label{r22k}
&&P^x((\wedge_{s,t}(A))^c,\ t<\zeta)=P^x\left(\left(\bigcap_{n=1}^{\infty}\wedge_{s,t}(A_n)\right)^c,\ t<\zeta\right),\nonumber\\
 &&Q^x((\wedge_{s,t}(A))^c,\ t<\zeta)=Q^x\left(\left(\bigcap_{n=1}^{\infty}\wedge_{s,t}(A_n)\right)^c,\ t<\zeta\right).
\end{eqnarray}

If $P^x(T_A<t,\ t<\zeta)=0$, then $P^x(s+D_A\circ
\theta_s\leq t-\frac{1}{l},\ t<\zeta)=0$ for
$l\in \mathbb{N}$. By (\ref{import}), we get $P^x(\wedge_{s,t-\frac{1}{l}}(A),\ t<\zeta)=0$ for $l\in \mathbb{N}$ satisfying $\frac{1}{l}< t-s$. Note that $\bigcap_{n=1}^{\infty}\wedge_{s,t-\frac{1}{l}}(A_n)\in
\mathcal{F}^0_{t}$. By (\ref{r22}) and the assumption  $Q^x|_{\cal{F}_t^0}$ is absolutely
continuous with respect to $P^x|_{\cal{F}_t^0}$ on $\{t<\zeta\}$, we have that
$Q^x(\wedge_{s,t-\frac{1}{l}}(A),\ t<\zeta)$=0. Then, we obtain by  (\ref{import}) that
$Q^x(s+D_A\circ \theta_s< t-\frac{1}{l},\ t<\zeta)=0$. Letting $l\rightarrow\infty$, we get $Q^x(s+D_A\circ \theta_s< t,\ t<\zeta)=0$.
Since $s\in (0,t)$ is arbitrary, (i) holds
by (\ref{lem-3.3-b-0}).

We now prove that (ii) holds. It is sufficient to show that for
any $t>0$, $P^x(T_A\ge t,\ t<\zeta)=0$ implies that $Q^x(T_A> t,\ t<\zeta)=0$ since it holds that
$$
\bigcup_{n=1}^{\infty}\left\{T_A\geq t+\frac{1}{n},\ t+\frac{1}{n}<\zeta\right\}=\{T_A>t,\ t<\zeta\}=\bigcup_{n=1}^{\infty}\left\{T_A> t+\frac{1}{n},\ t+\frac{1}{n}<\zeta\right\}.
$$

Suppose that $P^x(T_A\ge t,\ t<\zeta)=0$. Then, we obtain by
(\ref{lem-3.3-g-1}) that
\begin{eqnarray*}
 P^x\left(\frac{1}{l}+D_A\circ \theta_{\frac{1}{l}}\geq
t,\ t<\zeta\right)\downarrow0\ \ {\rm as}\ l\uparrow\infty,
\end{eqnarray*}
which together with (\ref{import2}) implies that
\begin{eqnarray}\label{r77}
P^x\left(\left(\wedge_{\frac{1}{l},t}(A)\right)^c,\ t<\zeta\right)\downarrow0\
\ {\rm as}\ l\uparrow\infty.
\end{eqnarray}
Note that
$\left(\bigcap_{n=1}^{\infty}\wedge_{\frac{1}{l},t}(A_n)\right)^c\in
\mathcal{F}^0_{t}$ for $l>\frac{1}{t}$. Then, we obtain by (\ref{import2}), (\ref{r22k}),
(\ref{r77}) and the assumption $Q^x|_{\cal{F}_t^0}$ is absolutely
continuous with respect to $P^x|_{\cal{F}_t^0}$ on $\{t<\zeta\}$ that
\begin{eqnarray}\label{r1000}
\lim_{l\to\infty} Q^x\left(\frac{1}{l}+D_A\circ
\theta_{\frac{1}{l}}> t,\ t<\zeta\right)\le\lim_{l\to\infty}
Q^x\left(\left(\wedge_{\frac{1}{l},t}(A)\right)^c,\ t<\zeta\right)=0.
\end{eqnarray}
Therefore, we obtain by (\ref{lem-3.3-g-1}) and (\ref{r1000}) that
$$
Q^x(T_A\ge t+\varepsilon,\ t<\zeta)=0,\ \ \forall \varepsilon>0,
$$
which implies that $Q^x(T_A> t,\ t<\zeta)=0$.
\hfill\fbox

\subsection{Proof of Theorem \ref{thm-2}}

(i)\ \ By
Lemma \ref{rlem-3.5}, it is sufficient to show that any thin Borel
measurable set for $Y$ is also polar for $Y$. Let $A\in {\cal B}$
be a thin set for $Y$. Define
 $$S=\{\omega:   T_A(\omega)=0\}
 $$
 and
 \begin{equation}\label{RL}
 R=\{\omega: T_A(\omega)<\infty\}.
 \end{equation}
 Then, $Q^x(S)=0$ for any $x\in E$.

By \cite[Theorem I.10.7]{BG68}, we know that $T_A$
is a stopping time relative to $\{\mathcal{F}_t^{X}\}$. Then, we
obtain by the Blumenthal's 0-1 law (\cite[Theorem I.5.17]{BG68})
that  $P^x(S)=0$ or 1 for any $x\in E$.

If  $P^x(S)=1$, then
 $P^x(S^c)=0$. Note that
 \begin{eqnarray*}
 S^c=\bigcup_{n=1}^{\infty}\left\{T_A>  \frac{1}{n}\right\}.
 \end{eqnarray*}
 Then, $P^x(T_A> \frac{1}{n})=0$ for $n\in \mathbb{N}$. By Lemma \ref{lem-3.3}, we get $Q^x(T_A> \frac{1}{n},\ \frac{1}{n}<\zeta)=0$ for $n\in \mathbb{N}$.
Since $Q^x(X_0=x)=1$, we obtain by the right continuity of $X$ that $ Q^x(\zeta>0)=1$. Thus, we have $$
 Q^x(S^c)=Q^x(T_A>0)=\lim_{n\to\infty}Q^x\left(T_A> \frac{1}{n},\ \frac{1}{n}<\zeta\right)=0,
 $$
 which implies that $Q^x(S)=1$. This contradicts with $Q^x(S)=0$. Therefore, $P^x(S)=0$ for any $x\in E$, which implies that $A$ is a thin set for $X$.

By the assumption  $X$ satisfies (H), we know that  $P^x(R)=0$ for
any $x\in E$.   For $t>0$, define
\begin{equation}\label{RR}
R_t=\{\omega: T_A(\omega)< t\}.
\end{equation}
Then, $P^x(R_t)=0$. By Lemma \ref{lem-3.3}, we obtain that $Q^x(R_t,\ t<\zeta)=0$ for $t>0$. Hence
$$
Q^x(R)\le\sum_{t\in\mathbb{Q}\cap (0,\infty)}Q^x(R_t,\ t<\zeta)=0.
$$
Since $x\in E$ is arbitrary, $A$ is a polar set for
$Y$.

(ii)\ \ By
Lemma \ref{rlem-3.5}, it is sufficient to show that any thin Borel
measurable set for $Y$ is also $m$-essentially polar for $Y$. Let
$A\in {\cal B}$ be a thin set for $Y$. Similar to (i), we obtain by Lemma \ref{lem-3.3} that $A$ is
a thin set for $X$.

Define $R$ and $R_t$, $t>0$, as in (\ref{RL}) and (\ref{RR}), respectively. By the assumption that $X$ satisfies ($H_m$), we get $P^m(R)=0$.
Then, $P^m(R_t)=0$ for $t>0$. By Lemma \ref{lem-3.3}, we
obtain that $Q^m(R_t,\ t<\zeta)=0$ for $t>0$. Hence
$$
Q^m(R)\le\sum_{t\in\mathbb{Q}\cap (0,\infty)}Q^m(R_t,\ t<\zeta)=0.
$$
\hfill\fbox

\section{Proof of Theorem \ref{thm-1}}\label{sec2}\setcounter{equation}{0}

The proof of (ii) is similar to that of (i). So we only
prove (i) below.

Let $X=(\Omega,X_t,P^x)$ be a standard process on $E$ and
$\tau=(\Theta,\tau_t,Q^0)$ be a subordinator with drift
coefficient $d>0$. We define $Y_t=X_{\tau_t}$ for $t\geq 0$. Then,
$Y=(\Omega\times\Theta,Y_t, P^x\times Q^0)$ is a standard process
on $E$.

$``\Rightarrow"$: Suppose $X$ satisfies (H). We will prove that $Y$ also satisfies (H).

(i) Suppose that $A\in {\cal B}$ is a polar set for $X$. Define
$$
\Omega_A=\{\omega\in\Omega: T_A<\infty\}.
$$
Then, for any $x\in E$,
\begin{eqnarray}\label{proof-thm-1.1-a}
P^x(\Omega_A)=0.
\end{eqnarray}
Since $d>0$, we have that
\begin{eqnarray}\label{proof-thm-1.1-b}
&&\ \ \ \{(\omega,w)\in \Omega\times \Theta: \exists\, t>0\ \ s.t.\ \ Y_t(\omega,w)\in A\}\nonumber\\
&&=\{(\omega,w)\in \Omega\times \Theta: \exists\, t>0\ \ s.t.\ \ X_{\tau_t(w)}(\omega)\in A\}\nonumber\\
&&\subset \Omega_A\times \Theta.
\end{eqnarray}
By (\ref{proof-thm-1.1-a}), (\ref{proof-thm-1.1-b}) and the
independence of $X$ and $\tau$, we find that for any $x\in E$,
$$
P^x\times Q^0(\{(\omega,w)\in \Omega\times \Theta: \exists\, t>0\ \ s.t.\ \ Y_t(\omega,w)\in A\})=0.
$$
Hence $A$ is also a polar set for $Y$.

(ii) For $y>0$, define $\varsigma(y)=\inf\{t\geq 0: \tau_t>y\}$. By Bertoin \cite[Theorem III.5]{B96}, we have that
\begin{eqnarray}\label{proof-thm-1.1-f}
Q^0(\tau_{\varsigma(y)}=y)=du(y),\ \ \forall y>0,
\end{eqnarray}
$u$ is continuous and positive on $(0,\infty)$, and
\begin{eqnarray}\label{r1}
u(0+)=1/d.
\end{eqnarray}

Let $A\in {\cal B}$. Suppose that $\omega\in\Omega$ satisfying
$T_A(\omega)=0$. By (\ref{proof-thm-1.1-f}) and (\ref{r1}), we get
$$
Q^0(\{w: \exists\, s>0\ \ s.t.\ \ \tau_s(w)\in \{0< t\le \epsilon:X_t(\omega)\in A\}\})=1,\ \ \forall \epsilon>0.
$$
Further, since $\tau_s(w)\ge ds$, we get
\begin{eqnarray}\label{r2}
Q^0(\{w: \exists\, 0<s\le \epsilon\ \ s.t.\ \ \tau_s(w)\in \{t:X_t(\omega)\in A\}\})=1,\ \ \forall \epsilon>0.
\end{eqnarray}

Suppose that $A\in {\cal B}$ is  not a polar set for $X$. Since
$X$ satisfies (H), $A$ is not a thin set for $X$. Hence $A$ has at
least one regular point, i.e., there exists some $x\in E$ such
that
\begin{eqnarray}\label{proof-thm-1.1-c}
P^x(T_A=0)=1.
\end{eqnarray}
By (\ref{r2}), (\ref{proof-thm-1.1-c}) and the independence of $X$
and $\tau$, we get
$$
 P^x\times Q^0(\{(\omega,w): T_A(\omega,w)=0\})=1.
 $$
Then, $x$ is a regular point of $A$ with respect to $Y$ and hence $A$ is not a thin set for $Y$.

We now show that $Y$ satisfies (H). By Lemma \ref{rlem-3.5}, we
need only show that any thin set  for $Y$ is also polar for $Y$.
Suppose that $A\in {\cal B}$ is a thin set for $Y$. If $A$ is not
polar for $Y$, then we obtain by (i) that $A$ is also not polar
for $X$. Further, we obtain by (ii) that $A$ is not a thin set for
$Y$. We have arrived at a contradiction.

$``\Leftarrow"$: Suppose that $Y$ satisfies (H). We will prove
that $X$ also satisfies (H).

If $X$ does not satisfy (H), then we obtain by Lemma
\ref{rlem-3.5} that there exits a $B\in {\cal B}$ such that $B$ is
a thin set   for $X$ but not a polar set for $X$. Since $d>0$, we
obtain by the independence of $X$ and $\tau$ that $B$ is a thin
set  for $Y$. By the assumption that $Y$ satisfies (H), we
conclude that $B$ is a polar set for $Y$.

Suppose that $\omega\in\Omega$ satisfying $T_B(\omega)<\infty$. By
(\ref{proof-thm-1.1-f}) and the fact that $u$ is positive on
$(0,\infty)$, we get
\begin{eqnarray}\label{r3}
Q^0(\{w: \exists\, s>0\ \ s.t.\ \ \tau_s(w)\in \{t>0:X_t(\omega)\in B\}\})>0.
\end{eqnarray}
Since $B$ is  not a polar set for $X$, there exists $x\in E$ such that $P^x(T_B(\omega)<\infty)>0$. Then, we obtain by (\ref{r3}) and the independence of $X$ and $\tau$ that
$$
 P^x\times Q^0(\{(\omega,w): T_B(\omega,w)<\infty\})>0,
 $$
which implies that $B$ is not a polar set for $Y$. We have arrived
at a contradiction.\hfill\fbox

\section{Invariance of (H) under absolutely continuous measure change}\label{sec4}\setcounter{equation}{0}

In this section, we apply Theorem \ref{thm-2} to study Hunt's
hypothesis (H) for standard processes. By virtue of absolutely
continuous measure change, we will give new examples of
subprocesses, L\'evy processes, and jump-diffusion processes
satisfying (H).

There exists a vast literature on the absolute continuity of
Markov processes. Skorohod \cite{Sk57,Sk65},  Kunita and Watanabe
\cite{KW67},  Newman \cite{Ne72,Ne73} and Jacod and Shiryaev
\cite{JS03} characterize the absolute continuity for L\'evy
processes. It\^{o} and Watanabe \cite{IW65}, Kunita
\cite{Ku69,Ku76} and Palmowski and Rolski \cite{PR02} discuss the
absolute continuity for general Markov processes. Dawson
\cite{Da68}, Liptser and Shiryaev \cite{LS77} and Kabanov, Liptser
and Shiryaev \cite{KLS80} study absolute continuity of solutions
to stochastic differential equations. We refer the reader to
Cheridito, Filipovi\'{c} and Yor \cite{CFY05} for a nice summary
of the references.

\subsection{Subprocesses (killing transformation)}

In this subsection, we will show that a standard process satisfies
(H) implies  that any of its standard subprocess satisfies (H).
First, let us recall some definitions from Blumenthal and Getoor
\cite{BG68}. Consider the following objects:
\begin{itemize}
\item[(i)] $W$: the space of all maps $w: [0,\infty]\to E_{\Delta}$
such that $w(\infty)=\Delta$ and if $w(t)=\Delta$ then
$w(s)=\Delta$ for all $s\geq t$.

\item[(ii)] Let $C_t$, $t\in [0,\infty]$, be the coordinate maps
$C_t(w)=w(t)$, and define in $W$ the $\sigma$-algebras
$\mathcal{C}^0=\sigma(C_t: t\in [0,\infty])$ and
$\mathcal{C}^0_t=\sigma(C_s: s\leq t)$ for $0\le t<\infty$.

\item[(iii)] Let $\varphi_t: W\to W$ be defined by
$\varphi_tw(s)=w(t+s)$.
\end{itemize}
Denote by $b{\cal B}$ and $b{\cal B}^*$ the sets of all bounded measurable and bounded universally measurable functions on $(E,\cal{B})$, respectively.

\begin{defi}\label{defi-3.2} (\cite{BG68}) (i) A Markov process $X=(\Omega,\mathcal{M},\mathcal{M}_t, X_t,\theta_t,P^x)$
with state space $(E,\mathcal{B})$ is said to be of {function
space type} provided $\Omega=W$, $\mathcal{M}\supset
\mathcal{C}^0$, $\mathcal{M}_t\supset \mathcal{C}_t^0$, $X_t=C_t$,
and $\theta_t=\varphi_t$.

(ii) Two Markov processes with the same state space $(E,{\cal B})$ are equivalent if they have the same transition function.

(iii) Let $X$ and $Y$ be two Markov processes with state space
$(E,{\cal B})$. Denote by $(P_t)$ and $(Q_t)$ the  transition
semigroups of $X$ and $Y$, respectively. $Y$ is called a
subprocess of $X$ if $Q_tf(x)\le P_tf(x)$ for all $x\in E$, $t\ge
0$ and $f\ge 0$ in $b{\cal B}^*$.
\end{defi}

\begin{lem}\label{equi}
Any standard process $X$ with state space $(E,{\cal B})$ is
equivalent to a standard process $C$ of function space type with
state space $(E,\mathcal{B})$. Moreover, $X$ satisfies (H) if and
only if $C$ satisfies (H).
\end{lem}
{\bf Proof.} Let $X=(\Omega,\mathcal{M},\mathcal{M}_t,
X_t,\theta_t,P^x)$ be a standard process with state space
$(E,{\cal B})$. Using the notation developed above Definition
\ref{defi-3.2} we define a map $\pi:\Omega\rightarrow W$ by
$(\pi\omega)(t)=X_t(\omega)$. Then, $C_t\circ\pi=X_t$. We define
measures $\hat{P}^x$ on $(W,\mathcal{C}^0)$ by
$\hat{P}^x=P^x\pi^{-1}$. By \cite[Theorem I.4.3]{BG68}, we know
that $C=(W,\mathcal{C}^0,\mathcal{C}^0_t,
C_t,\varphi_t,\hat{P}^x)$ is a Markov process equivalent to $X$.
Define $\mathcal{C}$ to be the completion of $\mathcal{C}^0$  with respect
to $\{\hat{P}^{\mu}: \mu\ \mbox{is a finite measure on}\ E_{\Delta}\}$.
For $t\in [0,\infty)$, define $\mathcal{C}_t$ to be the
completion of $\mathcal{C}^0_t$ in $\mathcal{C}$ with respect
to $\{\hat{P}^{\mu}: \mu\ \mbox{is a finite measure on}\
E_{\Delta}\}$. By the assumption
$X=(\Omega,\mathcal{M},\mathcal{M}_t, X_t,\theta_t,P^x)$ is a
standard process, we can check that
$C=(W,\mathcal{C},\mathcal{C}_t, C_t,\varphi_t,\hat{P}^x)$ is also
a standard process. Further, we can show that $X$ satisfies (H) if
and only if $C$ satisfies (H) by virtue of Lemma \ref{rlem-3.5}.
\hfill\fbox

\begin{thm}\label{dav}
Let $X$ be a standard process with state space $(E,{\cal B})$. If $X$ satisfies (H), then any standard subprocess of $X$ satisfies (H).
\end{thm}
{\bf Proof.} Suppose that $X=(X_t,P^x)$ satisfies (H) and
$Y=(Y_t,Q^x)$ is a standard subprocess of $X$.  By Lemma
\ref{equi}, we know that $X=(X_t,P^x)$ and $Y=(Y_t,Q^x)$ are
equivalent to some standard processes $C^X=(C_t,\hat{P}^x)$ and
$C^Y=({C}_t,\hat{Q}^x)$ of function space type with state space
$(E,\mathcal{B})$, respectively. Moreover, $C^X$ satisfies (H).
Since $\hat{Q}^x(A)\le \hat{P}^x(A)$ for any $A\in {\cal C}^0$ and
$x\in E$, $C^Y$ is absolutely continuous with respect to $C^X$.
Then, we obtain by Theorem \ref{thm-2} that ${C}^Y$ satisfies (H).
Therefore, $Y$ satisfies (H) by Lemma \ref{equi}.\hfill\fbox

\begin{rem}
Let $X=(X_t,P^x)$ be a standard process with state space
$(E,\cal{B})$. Suppose that $X$ satisfies (H). Let  $M=(M_t,0\leq
t<\infty)$ be a right continuous multiplicative functional (MF) of
$X$ satisfying $M_0=1$, $0\le M_t\le 1$ and $M_t\in {\cal F}^0_t$
for all $t$. By \cite[Corollary III.3.16]{BG68}, there is a
standard subprocess $\hat{X}=(\hat{X}_t,\hat{P}^x)$  with state
space $(E,\cal{E})$ such that
$\hat{E}^x[f(\hat{X}_t)]=E^x[f(X_t)M_t]$ for any $f\in b{\cal
B}^*$. By Theorem \ref{dav}, $\hat{X}$ satisfies (H). A concrete
example of the MF $M=(M_t,0\leq t<\infty)$ is given by
$$
M_t=\exp\left(-\int_0^tg(X_s)ds\right),\ \ t\geq 0,
$$
where $g\in b{\cal B}$ and $g\geq 0$.
\end{rem}

\subsection{L\'{e}vy processes (density transformation)}

Let $d\geq 1$ and $\mathbf{D}=D([0,\infty),\mathbb{R}^d)$ be the
space of mappings $\xi$ from $[0,\infty)$  into $\mathbb{R}^d$
which are right-continuous and have left limits. Denote
$x_t(\xi)=\xi(t)$. Define $\cal{F}_{\mathbf{D}}=\sigma(x_t,t\in
[0,\infty))$ and $\cal{F}^0_t=\sigma(x_s,s\in [0,t])$ for $t\in
[0,\infty)$. Any L\'{e}vy process on $\mathbb{R}^d$ induces a Hunt
process on $(\mathbf{D},\cal{F}_{\mathbf{D}})$ (cf. \cite[Section
40]{Sa99}). By Lemma \ref{equi}, when we consider (H) for a
L\'evy process, we may assume without loss of generality that it
is of the form $(x_t,P^x)$, where $P^x$ is a probability measure
defined on $(\mathbf{D},\cal{F}_{\mathbf{D}})$ for $x\in
\mathbb{R}^d$. Denote by $m_d$ the Lebesgue measure on $\mathbb{R}^d$.

In the sequel, we use  $\Phi$ or $(a, Q, \mu)$ to
denote the L\'{e}vy-Khintchine exponent of a L\'evy process $(X_t,P^x)$, which means that
\begin{eqnarray*}
E^0[\exp\{i\langle z,X_t\rangle\}]=\exp\{-t\Phi(z)\},\  z\in
\mathbb{R}^d,\ t\ge 0,
\end{eqnarray*}
where $E^0$ denotes the expectation {with respect to} $P^0$, and
\begin{eqnarray*}
\Phi(z)=i\langle a,z\rangle+\frac{1}{2}\langle z,Qz\rangle+\int_{\mathbb{R}^d}
\left(1-e^{i\langle z,x\rangle}+i\langle z,x\rangle 1_{\{|x|<1\}}\right)\mu(dx).
\end{eqnarray*}
Hereafter $\langle\cdot,\cdot\rangle$ and $|\cdot|$ denote the
standard Euclidean inner product and norm on $\mathbb{R}^d$,
respectively.

\begin{thm}\label{thm-4.2} (\cite[Chapter IV, Theorem 4.39 c]{JS03})
Let $(x_t, P^x)$ and $(x_t, P'^{x})$ be two L\'{e}vy processes on
$\mathbb{R}^d$ with L\'{e}vy-Khintchine exponents   $(a,Q,\mu)$
and $(a',Q',\mu')$, respectively. Then the following two
conditions are equivalent.

(1) $ P'^{0}|_{\cal{F}^0_t}\ll P^0|_{\cal{F}^0_t}$ for every $t\in (0,\infty)$.

(2) $Q=Q', \mu'\ll \mu$ with the function $K(x)$ defined by
$\frac{d\mu'}{d\mu}=K(x)$ satisfying
\begin{eqnarray}\label{thm-4.2-a}
\int_{\mathbb{R}^d}\left(1-\sqrt{K(x)}\right)^2\mu(dx)<\infty,
\end{eqnarray}
and
\begin{eqnarray}\label{thm-4.2-b}
a'-a+\int_{\{|x|<1\}}x(\mu'-\mu)(dx)\in \cal{R}(Q),
\end{eqnarray}
where $\cal{R}(Q):=\{Qx: x\in \mathbb{R}^d\}$.
\end{thm}

\begin{rem}
(i) Note that finiteness of the integral appearing in (\ref{thm-4.2-b}) follows from (\ref{thm-4.2-a}) (cf. \cite[Remark 33.3]{Sa99}).

(ii) If we let $h(x)=1_{\{|x|<1\}}$ in \cite[Chapter IV, Theorem
4.39 c]{JS03}, then $b$ and $b'$ of \cite[Chapter IV, Theorem 4.39
c]{JS03} satisfy $b=-a$ and $b'=-a'$.
\end{rem}

Combining Theorems \ref{thm-2} and \ref{thm-4.2}, we obtain the
following theorem.

\begin{thm}\label{main1}
Let $X$ and $X'$ be two L\'evy processes on $\mathbb{R}^d$ with
L\'{e}vy-Khintchine exponents $(a,Q,\mu)$ and $(a',Q',\mu')$,
respectively.  Suppose that Condition (2) in Theorem \ref{thm-4.2}
holds. Then,

(i) $X$ satisfies (H) implies that $X'$ satisfies (H).

(ii) $X$ satisfies ($H_{m_d}$) implies that $X'$ satisfies ($H_{m_d}$).
\end{thm}

It is well-known that any symmetric L\'evy process $X$ satisfies
($H_{m_d}$) (cf. \cite[Theorem 4.1.3]{FOT94} and \cite{Si77}). As a direct consequence of Theorem \ref{main1}, we have the following result on (H).

\begin{cor}\label{cor1}
Let $Q$ be a symmetric nonnegative-definite $d\times d$ matrix and
$\mu$ be a measure on $\mathbb{R}^d$  satisfying $\mu(\{0\})=0$,
$\int_{\mathbb{R}^d}(1\wedge|x|^2)\mu(dx)<\infty$, and
$\mu(A)=\mu(-A)$ for any $A\in {\cal B}(\mathbb{R}^d)$. Suppose
that $K$ is a nonnegative measurable function on $\mathbb{R}^d$
satisfying
$$
\int_{\mathbb{R}^d}\left(1-\sqrt{K(x)}\right)^2\mu(dx)<\infty,
$$
and $a\in \mathbb{R}^d$ satisfying
$$
a+\int_{\{|x|< 1\}}x(K(x)-1)\mu(dx)\in \cal{R}(Q).
$$
Let $Y$ be a L\'evy process on  $\mathbb{R}^d$ with
L\'{e}vy-Khintchine exponent $(a,Q,K(x)d\mu)$. Then $Y$ satisfies
($H_{m_d}$). If in addition $Y$ has resolvent densities with
respect to $m_d$, then $Y$ satisfies (H).
\end{cor}

Now we give a useful lemma on the absolute continuity of measure change for L\'evy processes.

\begin{lem}\label{lem-4.2-1}
Let $(x_t, P^x)$ be a L\'{e}vy process  on $\mathbb{R}^d$ with
L\'{e}vy-Khintchine exponent  $(a,Q,\mu)$ and $\mu_1$ be a
measure on $\mathbb{R}^d\backslash\{0\}$.

(i) Suppose $\mu_1\leq \mu$ with the function $k(x)$ defined by
$k(x)=\frac{d\mu_1}{d\mu}$ satisfying $\int_{\{|x|<
1\}}k^2(x)\mu(dx)<\infty$. Denote $\mu':=\mu-\mu_1$ and let $(x_t,
P'^{x})$ be a L\'{e}vy process on $\mathbb{R}^d$ with
L\'{e}vy-Khintchine exponent $(a',Q,\mu')$, where
$a':=a+\int_{\{|x|< 1\}}x\mu_1(dx)$. Then $
P'^{x}|_{\cal{F}^0_t}\ll P^x|_{\cal{F}^0_t}$ for $x\in
\mathbb{R}^d$ and $t>0$.

(ii) Suppose $\mu_1 \ll\mu$ with the function $k(x)$ defined by
$k(x)=\frac{d\mu_1}{d\mu}$ satisfying $\int_{\{|x|<
1\}}k^2(x)\mu(dx)<\infty$ and $\int_{\{|x|\geq
1\}}k(x)\mu(dx)<\infty$. Denote $\mu'':=\mu+\mu_1$ and let $(x_t,
P''^{x})$ be a L\'{e}vy process on $\mathbb{R}^d$ with
L\'{e}vy-Khintchine exponent $(a'',Q,\mu'')$, where
$a'':=a-\int_{\{|x|< 1\}}x\mu_1(dx)$. Then $
P''^{x}|_{\cal{F}^0_t}\ll P^x|_{\cal{F}^0_t}$ for $x\in
\mathbb{R}^d$ and $t>0$.
\end{lem}
{\bf Proof.} Note that
$$\int_{\{|x|< 1\}}|x|\mu_1(dx)\le\left(\int_{\{|x|< 1\}}|x|^2\mu(dx)\right)^{1/2}\left(\int_{\{|x|< 1\}}k^2(x)\mu(dx)\right)^{1/2}<\infty.
$$
Obviously, condition (\ref{thm-4.2-b}) holds. By Theorem \ref{thm-4.2},  it is sufficient to show that condition (\ref{thm-4.2-a}) holds.

(i) Denote $K(x)=\frac{d\mu'}{d\mu}$. Then, $K(x)=1-k(x)$. By the
assumption that $\mu_1\leq \mu$, we get $0\leq k(x)\leq 1$.  It
follows that $\sqrt{1-k(x)}\geq 1-k(x)$, which implies that
$(1-\sqrt{1-k(x)})^2\leq k^2(x)$. Therefore, $\int_{\{|x|<
1\}}k^2(x)\mu(dx)<\infty$ implies that $
\int_{\mathbb{R}^d}\left(1-\sqrt{K(x)}\right)^2\mu(dx)<\infty. $

(ii) Denote $K(x)=\frac{d\mu'}{d\mu}$. Then, $K(x)=1+k(x)$, where
$k(x)\geq 0$.  It follows that $\sqrt{1+k(x)}\leq  1+k(x)$,  which
implies that $(\sqrt{1+k(x)}-1)^2\leq k^2(x)$. Therefore,
$\int_{\{|x|< 1\}}k^2(x)\mu(dx)<\infty$ implies that $
\int_{\mathbb{R}^d}\left(1-\sqrt{K(x)}\right)^2\mu(dx)<\infty. $
\hfill\fbox

\begin{cor}\label{cor-4.2-2}
Let $(x_t, P^x)$ be a L\'{e}vy process  on $\mathbb{R}^d$ with
L\'{e}vy-Khintchine exponent  $(a,Q,\mu)$. Suppose  $\mu_1$ is a
finite measure on $\mathbb{R}^d\backslash\{0\}$ such that
$\mu_1\leq \mu$. Denote $\mu':=\mu-\mu_1$ and let $(x_t, P'^{x})$
be a L\'{e}vy process on $\mathbb{R}^d$ with L\'{e}vy-Khintchine
exponent $(a',Q,\mu')$, where $a':=a+\int_{\{|x|< 1\}}x\mu_1(dx)$.
Then $ P'^{x}|_{\cal{F}^0_t}\ll P^x|_{\cal{F}^0_t}$ for $x\in
\mathbb{R}^d$ and $t>0$.
\end{cor}
{\bf Proof.} Denote $k(x)=\frac{d\mu_1}{d\mu}$. Then, $0\leq k(x)\leq 1$, which together with the fact that $\mu_1$ is a finite measure implies that
\begin{eqnarray*}
\int_{\{|x|< 1\}}k^2(x)\mu(dx)&\leq& \int_{\{|x|< 1\}}k(x)\mu(dx)=\int_{\{|x|< 1\}}\mu_1(dx)<\infty.
\end{eqnarray*}
Thus, the proof is complete by Lemma \ref{lem-4.2-1}(i).\hfill\fbox

Combining \cite[Theorem 3.1 and Proposition 3.3]{HS18}, Theorem
\ref{thm-2} and Corollary \ref{cor-4.2-2}, we obtain the following
result, which implies that big jumps have no effect on the
validity of (H) for any L\'{e}vy process.

\begin{pro}\label{thm-4.7}
Let $X$ be a L\'{e}vy process  on $\mathbb{R}^d$ with
L\'{e}vy-Khintchine exponent  $(a,Q,\mu)$. Suppose  that $\mu_1$
is a finite measure on $\mathbb{R}^d\backslash\{0\}$ such that
$\mu_1\leq \mu$. Denote $\mu':=\mu-\mu_1$ and let $X'$ be a
L\'{e}vy process on $\mathbb{R}^d$ with L\'{e}vy-Khintchine
exponent $(a',Q,\mu')$, where $a':=a+\int_{\{|x|< 1\}}x\mu_1(dx).$
Then,

(i) $X$ and $X'$ have same semipolar sets.

(ii) $X$ and $X'$ have same $m_d$-essentially polar sets.

(iii) $X$ satisfies (H) if and only if $X'$ satisfies (H).

(iv) $X$ satisfies $(H_{m_d})$ if and only if $X'$ satisfies $(H_{m_d})$.
\end{pro}

Finally, we give a new class of purely discontinuous L\'evy processes satisfying (H).
\begin{pro}
Suppose that $\rho_1$ is a nonnegative measurable
function on $\mathbb{R}^d$ satisfying $\rho_1(x)=\rho_1(-x)$ for $x\in \mathbb{R}^d$ and
$$\int_{\mathbb{R}^d}(1\wedge|x|^2)\rho_1(x)dx<\infty,$$ and $\rho_2$ is a
nonnegative measurable function on $\mathbb{R}^d$ satisfying
$$\int_{\{|x|< 1\}}\rho_2^2(x)\rho_1(x)dx<\infty\ \ {\rm and}\ \ \int_{\{|x|\geq  1\}}\rho_2(x)\rho_1(x)dx<\infty. $$
Denote $a=-\int_{\{|x|< 1\}}x\rho_2(x)\rho_1(x)dx$ and
$\mu=(1+\rho_2(x))\rho_1(x)dx$. Let $Y$ be a L\'evy process on
$\mathbb{R}^d$ with L\'{e}vy-Khintchine exponent $(a,0,\mu)$. Then
$Y$ satisfies (H).
\end{pro}
{\bf Proof.} Let $X$ be a L\'evy
process on $\mathbb{R}^d$ with L\'{e}vy-Khintchine exponent
$(0,0,\rho_1(x)dx)$.  If $\rho_1(x)dx$ is a finite measure, then $X$ is a
compound Poisson process and hence satisfies (H). We now assume
that $\rho_1(x)dx$ is an infinite measure.  Then, $X$  is a symmetric L\'evy process
which has transition densities by \cite[Theorem 27.7]{Sa99}. So
$X$ satisfies (H). Therefore, the proof is complete by Theorem
\ref{thm-2} and Lemma \ref{lem-4.2-1}(ii).\hfill\fbox

\subsection{Jump-diffusion processes (supermartingale transformation)}

Throughout this subsection, we make the following assumptions.

\noindent (A.1) $a=(a_{ij})$ is a bounded continuous mapping on $\mathbb{R}^d$ with
values in the set of positive definite symmetric $d\times d$ matrices such that $\frac{\partial a_{ij}}{\partial x_i}$ is a bounded measurable function on $\mathbb{R}^d$ for $1\le i,j\le d$.

\noindent (A.2) $b$ is a bounded measurable mapping on $\mathbb{R}^d$
with values in $\mathbb{R}^d$.

\noindent (A.3)  $\gamma$ is a nonnegative symmetric
measurable function on $\mathbb{R}^d\times
\mathbb{R}^d\backslash{\rm diagonal}$ such that
$$
\int_{\Gamma}\frac{|y|^2}{1+|y|^2}\gamma(x,x+y)dy
$$
is bounded continuous for all $\Gamma\in{\cal B}(\mathbb{R}^d\backslash\{0\})$.

\noindent (A.4)   $k$ is a nonnegative measurable function on $\mathbb{R}^d\times
\mathbb{R}^d\backslash{\rm diagonal}$ such that
\begin{equation}\label{yor}
\int_{\mathbb{R}^d}k^2(x,y)\gamma(x,y)dy\ {\rm is\ bounded\ on}\ \mathbb{R}^d,
\end{equation}
and
$$
\int_{\Gamma}\frac{|y|^2}{1+|y|^2}k(x,x+y)\gamma(x,x+y)dy
$$
is bounded continuous for all $\Gamma\in{\cal B}(\mathbb{R}^d\backslash\{0\})$.

We use $C^{\infty}_c(\mathbb{R}^d)$ to denote the space of smooth
functions on $\mathbb{R}^d$ with compact support. Define
\begin{eqnarray}\label{4.3-b-1}
\vartheta(x,y)=[1+k(x,y)]\gamma(x,y).
\end{eqnarray}
For $f\in C^{\infty}_c(\mathbb{R}^d)$, define
\begin{eqnarray*}
\mathcal{A}^{a,{b},\vartheta}f(x)&:=&\frac{1}{2}\sum_{i,j=1}^da_{ij}(x)\frac{\partial^2f(x)}{\partial x_i\partial x_j}+\sum_{i=1}^d{b}_i(x)\frac{\partial f(x)}{\partial x_i}\\
&&+\int_{\mathbb{R}^d}\left(f(x+y)-f(x)-\frac{\langle y,\nabla f(x)\rangle}{1+|y|^2}\right)\vartheta(x,x+y)dy.
\end{eqnarray*}
By (A.1)--(A.4) and Stroock \cite[Theorem 4.3 and Remark, page 232]{St75}, we know that the martingale problem for
$\mathcal{A}^{a,{b},\vartheta}$ is well-posed and the
corresponding solution is a strong Feller process. Denote by $(x_t,P^x)$ the Markov process
associated with $\mathcal{A}^{a,{b},\vartheta}$ on
$(\mathbf{D},\cal{F}_{\mathbf{D}})$, where $(\mathbf{D},\cal{F}_{\mathbf{D}})$ is the Skorohod space as defined in \S 4.2.

\begin{thm}\label{stroock}
Assume that (A.1)--(A.4) hold. Then $(x_t,P^x)$ satisfies (H).
\end{thm}
{\bf Proof.} For $f\in C^{\infty}_c(\mathbb{R}^d)$ and $x\in \mathbb{R}^d$, we define
\begin{eqnarray*}
\mathcal{A}^{a,\gamma}f(x)&:=&\frac{1}{2}\sum_{i,j=1}^d\frac{\partial}{\partial x_i}\left(a_{ij}\frac{\partial f}{\partial
x_j}\right)(x)+\int_{\mathbb{R}^d}\left(f(x+y)-f(x)-\frac{\langle y,\nabla
f(x)\rangle}{1+|y|^2}\right)\gamma(x,x+y)dy.
\end{eqnarray*}
Define
$$
\varsigma_i(x)=\sum_{j=1}^d\frac{\partial a_{ij}}{\partial x_j},\ \ 1\le i\le d,
$$
and set $\varsigma=(\varsigma_1,\dots,\varsigma_d)$. Then,
\begin{eqnarray*}
\mathcal{A}^{a,\gamma}f(x)&=&\frac{1}{2}\sum_{i,j=1}^da_{ij}(x)\frac{\partial^2f(x)}{\partial
x_i\partial
x_j}+\sum_{i=1}^d\varsigma_i(x)\frac{\partial f(x)}{\partial x_i}\nonumber\\
&&+\int_{\mathbb{R}^d}\left(f(x+y)-f(x)-\frac{\langle y,\nabla
f(x)\rangle}{1+|y|^2}\right)\gamma(x,x+y)dy.
\end{eqnarray*}
Hence we obtain by (A.1)--(A.3) and \cite[Theorem 4.3 and Remark, page 232]{St75} that the martingale problem for
$\mathcal{A}^{a,\gamma}$ is well-posed and the
corresponding solution is a strong Feller process. Denote by $(x_t,P^{*x})$ the Markov process
associated with $\mathcal{A}^{a,\gamma}$ on
$(\mathbf{D},\cal{F}_{\mathbf{D}})$.

We fix a bounded continuous function
$\chi:\mathbb{R}^d\rightarrow\mathbb{R}^d$ such that $\chi(y)=y$
on a neighborhood of $0$. For $x\in \mathbb{R}^d$, define
\begin{eqnarray}\label{4.3-d}
\beta(x)=\varsigma(x)+\int_{\mathbb{R}^d}\left(\chi(y)-\frac{y}{1+|y|^2}\right)\gamma(x,x+y)dy,
\end{eqnarray}
and
\begin{eqnarray}\label{4.3-f}
\tilde{\beta}(x)=b(x)+\int_{\mathbb{R}^d}\left(\chi(y)-\frac{y}{1+|y|^2}\right)
\vartheta(x,x+y)dy.
\end{eqnarray}
Then, for $f\in C^{\infty}_c(\mathbb{R}^d)$, we have that
\begin{eqnarray}\label{4.3-c}
\mathcal{A}^{a,\gamma}f(x)&=&\frac{1}{2}\sum_{i,j=1}^da_{ij}(x)\frac{\partial^2f(x)}{\partial
x_i\partial
x_j}+\sum_{i=1}^d\beta_i(x)\frac{\partial f(x)}{\partial x_i}\nonumber\\
&&+\int_{\mathbb{R}^d}\left(f(x+y)-f(x)-\langle \nabla
f(x),\chi(y)\rangle\right)\gamma(x,x+y)dy,
\end{eqnarray}
and
\begin{eqnarray}\label{4.3-e}
\mathcal{A}^{a,{b},\vartheta}f(x)&=&\frac{1}{2}\sum_{i,j=1}^da_{ij}(x)\frac{\partial^2f(x)}{\partial x_i\partial x_j}+\sum_{i=1}^d\tilde{\beta}_i(x)\frac{\partial f(x)}{\partial x_i}\nonumber\\
&&+\int_{\mathbb{R}^d}\left(f(x+y)-f(x)-\langle \nabla
(f),\chi(y)\rangle\right)\vartheta(x,x+y)dy.
\end{eqnarray}

Define
\begin{eqnarray}\label{4.3-b}
\phi(x)=a^{-1}(x)\left({b}(x)-\varsigma(x)-\int_{\mathbb{R}^d}\frac{y}{1+|y|^2}k(x,x+y)\gamma(x,x+y)dy\right),\
\ x\in \mathbb{R}^d.
\end{eqnarray}
Then, $\phi$ is locally bounded on $\mathbb{R}^d$. By
(\ref{4.3-b-1}), (\ref{4.3-d}), (\ref{4.3-f}) and (\ref{4.3-b}),
we get
\begin{eqnarray}\label{4.3-g}
\tilde{\beta}(x)&=&\beta(x)-\varsigma(x)+b(x)+\int_{\mathbb{R}^d}\left(\chi(y)-\frac{y}{1+|y|^2}\right)
\left(\vartheta(x,x+y)-\gamma(x,x+y)\right)dy\nonumber\\
&=&\beta(x)+\left[a(x)\phi(x)+\int_{\mathbb{R}^d}\frac{y}{1+|y|^2}k(x,x+y)\gamma(x,x+y)dy\right]\nonumber\\
&&+\int_{\mathbb{R}^d}\left(\chi(y)-\frac{y}{1+|y|^2}\right)k(x,x+y)\gamma(x,x+y)dy\nonumber\\
&=&\beta(x)+a(x)\phi(x)+\int_{\mathbb{R}^d}\chi(y)k(x,x+y)\gamma(x,x+y)dy.
\end{eqnarray}
Note that
(\ref{yor}) implies that
\begin{eqnarray}\label{4.3-h}
\int_{\mathbb{R}^d}[(k(x,y)+1)\log(k(x,y)+1)-k(x,y)]\gamma(x,x+y)dy\
\ {\rm is\ bounded\ on}\ \mathbb{R}^d.
\end{eqnarray}
Then, we obtain by \cite[Theorem 2.4 and Remark 2.5]{CFY05},
(\ref{4.3-b-1}), (\ref{4.3-c}), (\ref{4.3-e}), (\ref{4.3-g}) and
(\ref{4.3-h}) that $(x_t,P^x)$ is induced by a supermartingale
transformation of $(x_t,P^{*x})$. Hence $P^x|_{\cal{F}_t^0}$ is
absolutely continuous with respect to $P^{*x}|_{\cal{F}_t^0}$ for any
$x\in \mathbb{R}^d$ and $t>0$.

We now use the theory of Dirichlet forms to show that
$(x_t,P^{*x})$ satisfies (H).  The reader is referred to
Fukushima, Oshima and Takeda \cite{FOT94} for notation and
terminology used below. First, note that
$\mathcal{A}^{a,\gamma}f\in L^2(\mathbb{R}^d;dx)$ for $f\in
C^{\infty}_c(\mathbb{R}^d)$. In fact, suppose ${\rm
supp}[f]\subset \{x\in \mathbb{R}^d:|x|\leq N\}$ for some
$N\in\mathbb{N}$. Then, we obtain by (A.3) that \begin{eqnarray*}
&&\int_{\{|x|>2N\}}\left|\int_{\mathbb{R}^d}\left(f(x+y)-f(x)-\frac{\langle
y,\nabla
f(x)\rangle}{1+|y|^2}\right)\gamma(x,x+y)dy\right|^2dx\\
&=&\int_{\{|x|>2N\}}\left(\int_{\{|y|\ge1\}}{f(x+y)}\gamma(x,x+y)dy\right)^2dx\\
&\le&\left(\sup_{x\in \mathbb{R}^d}\int_{\{|y|\ge1\}}\gamma(x,x+y)dy\right)\int_{\{|x|>2N\}}\int_{\{|y|\ge1\}}{f^2(x+y)}\gamma(x,x+y)dydx\\
&\le&\left(\sup_{x\in \mathbb{R}^d}\int_{\{|y|\ge1\}}\gamma(x,x+y)dy\right)\int_{\mathbb{R}^d}\int_{\{|y-x|\ge1\}}{f^2(y)}\gamma(x,y)dydx\\
&=&\left(\sup_{x\in \mathbb{R}^d}\int_{\{|y|\ge1\}}\gamma(x,x+y)dy\right)\int_{\mathbb{R}^d}\int_{\{|y-x|\ge1\}}{f^2(y)}\gamma(x,y)dxdy\\
&\le&\left(\sup_{x\in \mathbb{R}^d}\int_{\{|y|\ge1\}}\gamma(x,x+y)dy\right)^2\int_{\mathbb{R}^d}f^2(y)dy\\
&<&\infty.
\end{eqnarray*}
We consider the symmetric bilinear form on $L^2(\mathbb{R}^d;dx)$:
\begin{eqnarray*}
{\cal
E}(f,g)&=&-\int_{\mathbb{R}^d}\mathcal{A}^{a,\gamma}f(x)g(x)dx\\
&=&\frac{1}{2}\sum_{i,j=1}^d\int_{\mathbb{R}^d}a_{ij}\frac{\partial f}{\partial
x_i}\frac{\partial g}{\partial
x_j}dx+\int_{\mathbb{R}^d}\int_{\mathbb{R}^d}(f(x)-f(y))(g(x)-g(y))dxdy,\ \ f,g\in C^{\infty}_c(\mathbb{R}^d).
\end{eqnarray*}
$({\cal E},C^{\infty}_c(\mathbb{R}^d))$ is closable on $L^2(\mathbb{R}^d;dx)$ and its closure $({\cal E},D({\cal E}))$ is a regular symmetric Dirichlet form.

Let
$(X^{\cal E}_t,P^{{\cal E},x})$ be a Hunt process
associated with $({\cal E},D({\cal E}))$. For
$f\in C^{\infty}_c(\mathbb{R}^d)$, define
$$
g=(1-\mathcal{A}^{a,\gamma})f.
$$
Then, $g$ is a bounded measurable function on $\mathbb{R}^d$ and $g\in L^2(\mathbb{R}^d;dx)$. Denote by $R^{\cal E}_1$ the 1-resolvent of $X^{\cal E}$ and define
$$
{\cal M}^f_t=R^{\cal E}_1g(X^{\cal E}_t)-R^{\cal E}_1g(X^{\cal
E}_0)-\int_0^t(R^{\cal E}_1g-g)(X^{\cal E}_s)ds, \ \ t\ge 0.
$$
It is known that $\{{\cal M}^{f}_t\}$ is a martingale under
$P^{{\cal E},x}$ for $x\in \mathbb{R}^d$ (cf. \cite[Chapter 4,
Proposition 1.7]{EK}). Denote by $G_1$ the 1-resolvent of ${\cal
E}$ and define
$$
M^{f}_t=f(X^{\cal E}_t)-f(X^{\cal
E}_0)-\int_0^t\mathcal{A}^{a,\gamma}f(X^{\cal E}_s)ds, \ \ t\ge 0.
$$
Since
$f=G_1g$ $dx$-a.e., we get $f=R^{\cal E}_1g$ q.e.. Hence $\{M^{f}_t\}$ is a martingale under $P^{{\cal
E},x}$ for q.e. $x\in \mathbb{R}^d$.

Let $\Psi$ be a countable subset of $C^{\infty}_c(\mathbb{R}^d)$
such that for any $f\in C^{\infty}_c(\mathbb{R}^d)$ there exist
$\{f_n\}\subset\Psi$ satisfying $\|f_n-f\|_{\infty},
\|\partial_if_n-\partial_if\|_{\infty},\|\partial_i\partial_jf_n-\partial_i\partial_jf\|_{\infty}\rightarrow0$
as $n\rightarrow\infty$ for any $i,j\in\{1,2,\dots,d\}$. Then, there is
an ${\cal E}$-exceptional set of $\mathbb{R}^d$, denoted by $F$,
such that $\{M^{f}_t\}$ is a martingale under $P^{{\cal E},x}$
for any $x\in F^c$.
We obtain by taking limits that $\{M^{f}_t\}$ is a martingale
under $P^{{\cal E},x}$ for any $f\in C^{\infty}_c(\mathbb{R}^d)$
and q.e. $x\in \mathbb{R}^d$. Thus, by the uniqueness of
solutions to the martingale problem for $\mathcal{A}^{a,\gamma}$, we conclude
that $(x_t,P^{*x})$ is a Hunt process associated with the symmetric Dirichlet form
$({\cal E},D({\cal E}))$. Since $(x_t,P^{*x})$ is a strong
Feller process,
$(x_t,P^{*x})$ must satisfy (H) (cf. \cite[Theorems 4.1.2 and
4.1.3]{FOT94}). Therefore, we obtain by Theorem \ref{thm-2} that
$(x_t,P^x)$ satisfies (H).\hfill\fbox

\section{Invariance of (H) under subordination and remark}\label{sec5}\setcounter{equation}{0}

As a direct application of Theorem \ref{thm-1}, we obtain the
following result.

\begin{pro}\label{cor-5.2}
Suppose that $X$ is a L\'{e}vy process on $\mathbb{R}^d$ with
L\'{e}vy-Khintchine exponent $\varphi$. Let $c>0$
be a constant and $\nu$ be a measure on $(0,\infty)$
satisfying $\int_0^{\infty}(1\wedge x)\nu(dx)<\infty$. Then $X$ satisfies (H) if and only if the
L\'{e}vy process on $\mathbb{R}^d$ with L\'{e}vy-Khintchine
exponent
\begin{eqnarray}\label{cor-5.2-a}
\Phi(z):=c\,\varphi(z)+\int_0^{\infty}\left(1-e^{-\varphi(z)x}\right)
\nu(dx),\ z\in \mathbb{R}^d
\end{eqnarray}
satisfies (H). In particular, for $0<\alpha<1$, $X$ satisfies (H) if and only if the L\'{e}vy
process on $\mathbb{R}^d$ with L\'{e}vy-Khintchine exponent
\begin{eqnarray*}
\Phi(z)=c\,\varphi(z)+(\varphi(z))^{\alpha},\ z\in \mathbb{R}^d
\end{eqnarray*}
satisfies (H).
\end{pro}
{\bf Proof.} Let  $\tau$ be a subordinator with drift coefficient
$c$ and L\'{e}vy measure $\nu$, which is independent of $X$.
Define $Y_t:=X_{\tau_t}$ for $t\ge 0$. Then, $Y=(Y_t)$ has the
L\'{e}vy-Khintchine exponent $\Phi$ defined by (\ref{cor-5.2-a}).
Therefore, $X$ satisfies (H) if and only if $Y$ satisfies (H) by Theorem \ref{thm-1}. The second
assertion is proved by letting $\tau_t=ct+\beta_t$, where $\beta$ is a stable subordinator of index $\alpha$ which is independent $X$.\hfill\fbox

Note that the uniform motion
$X_t=t$ on $\mathbb{R}$ does not satisfy (H). The sufficient part of Theorem \ref{thm-1} can be regarded as a generalization of  the following proposition.

\begin{pro}\label{cor-5.1}(cf. \cite[Proposition 1.6]{HS12})
Let $X$ be a subordinator. If $X$ satisfies (H), then its drift
coefficient equals 0.
\end{pro}

Proposition \ref{cor-5.1} can be extended to the $d$-dimensional
case as follows.

\begin{pro}\label{pro-5.6}
Let $X$ be a L\'{e}vy process on $\mathbb{R}^d$ with
L\'{e}vy-Khintchine exponent $(a,0,\mu)$ satisfying
$\int_{\mathbb{R}^d} (1\wedge|x|)\mu(dx)<\infty$. If $X$ satisfies
(H), then its drift coefficient equals 0.
\end{pro}
{\bf Proof.} The 1-dimensional case follows by \cite{Ke69,Br71}.
Now we consider the case  that $d>1$. Denote by $a'$ the drift
coefficient of $X$, i.e., $ a'=-\left(a+\int_{\{|x|<
1\}}x\mu(dx)\right). $ Then, we have
\begin{eqnarray}\label{lem-4.2-a}
E^0[e^{i\langle z,X_1\rangle}]&=&\exp\left[-\left(i\langle
-a',z\rangle+\int_{\mathbb{R}^d}\left(1-e^{i\langle
z,x\rangle}\right)\mu(dx)\right)\right],\ \ z\in \mathbb{R}^d.
\end{eqnarray}

Suppose that $X$ satisfies (H) and $a'\neq 0$. Define
$S=\{ta'|t\in \mathbb{R}\}$ and let $Y=(Y_t)$ be the projection
process of $X$ on $S$. Then, we obtain by \cite[Lemma 3.4]{HS18}
that  $Y$ is a one-dimensional L\'{e}vy process satisfying (H).

Denote by $P$ the projection operator from $\mathbb{R}^d$ to $S$.
Then, we obtain by (\ref{lem-4.2-a}) that for $z\in \mathbb{R}^d$,
\begin{eqnarray*}
E^0[e^{i\langle z,Y_1\rangle}]&=&E^0[e^{i\langle z,PX_1\rangle}]\\
&=&E^0[e^{i\langle Pz,X_1\rangle}]\nonumber\\
&=&\exp\left\{-\left(i\langle -a',Pz\rangle+\int_{\mathbb{R}^d}\left(1-e^{i\langle Pz,x\rangle}\right)\mu(dx)\right)\right\}\nonumber\\
&=&\exp\left\{-\left(i\langle -Pa',z\rangle+\int_{\mathbb{R}^d}\left(1-e^{i\langle z,Px\rangle}\right)\mu(dx)\right)\right\}\nonumber\\
&=&\exp\left\{-\left(i\langle
-a',z\rangle+\int_{\mathbb{R}^d}\left(1-e^{i\langle
z,y\rangle}\right)\mu_P(dy)\right)\right\},
\end{eqnarray*}
where $\mu_P$ is the image measure of $\mu$ under the map $P$.
Hence the drift coefficient of $Y$ is $a'$. Since the proposition
is true for the 1-dimensional case, $a'=0$. We have arrived at a
contradiction.\hfill\fbox

\begin{rem}\label{rem-5.3}
Let $X=(X_t)$ be a L\'{e}vy process  on $\mathbb{R}^d$ with
L\'{e}vy-Khintchine exponent $\Phi$ or $(a,Q,\mu)$. $X$ is called
a pure jump L\'{e}vy process if $Q=0$, $\int_{\mathbb{R}^d}(
1\wedge|x|)\mu(dx)<\infty$, and the drift coefficient
$a'=-\left(a+\int_{\{|x|< 1\}}x\mu(dx)\right)=0$. In this case,
$\Phi$ can be expressed by $ \Phi(z)=\int_{\mathbb{R}^d}
\left(1-e^{i\langle z,x\rangle}\right)\mu(dx)$ and $X$ can be
expressed by $X_t=\sum_{0<s\leq t}\triangle X_s,$ where $\triangle
X_s=X_s-X_{s-}$.

Based on Proposition \ref{pro-5.6}, it is natural to ask the
following question. \vskip 0.3cm \noindent {\bf Question} Does any
pure jump L\'{e}vy process satisfy (H)? \vskip 0.3cm \noindent If
the answer to the above question is affirmative for subordinators,
i.e., any pure jump subordinator satisfies (H), then Theorems
\ref{thm-3} and  \ref{thm-1} imply the following claim:

Let $(X_t)$ be a standard process on an LCCB state space and
$(\tau_t)$ be a subordinator which is independent of $(X_t)$. Then
$(X_{\tau_t})$ satisfies (H) if and only if either $(X_t)$ or
$(\tau_t)$ satisfies (H).
\end{rem}

\bigskip

{ \noindent {\bf\large Acknowledgments}  This work was supported
by  National Natural Science Foundation of China (Grant No.
11771309) and Natural Science and Engineering Research Council of
Canada.


\end{document}